\documentclass[11pt,a4paper]{article}

\usepackage{fancyhdr}
\usepackage{amsmath, amssymb}
\usepackage{fourier}
\usepackage[top=2.5cm,bottom=2.5cm,left=2.5cm,right=2.5cm]{geometry}
\usepackage{graphicx}
\usepackage[colorlinks=true,urlcolor=blue,pdfstartview=FitV, linkcolor=blue,citecolor=blue,bookmarksopen, bookmarksnumbered={true}]{hyperref}
\usepackage{enumerate}
\newcommand\tr{\mathrm{tr}}
\newtheorem{theorem}{Theorem}[section]
\newtheorem{definition}{Definition}[section]
\newtheorem{corollary}{Corollary}[section]
\newtheorem{example}{Example}[section]

\newtheorem{lemma}{Lemma}[section]

\def\Q{\widetilde{Q}}

\title{Existence  of positive solutions for   generalized  Lyapunov equations via a coupled fixed point theorem}
\author{Maher Berzig and Bessem Samet}

\begin{document}

\maketitle
\hrule
\begin{abstract}
We consider the generalized continuous-time Lyapunov equation:
$$
A^*XB + B^*XA =-Q,
$$
 where $Q$ is an $N\times N$  Hermitian positive definite matrix and $A,B$ are arbitrary $N\times N$ matrices. Under some conditions, using the coupled fixed point theorem of Bhaskar and Lakshmikantham, we establish the existence and uniqueness of Hermitian positive definite solution for such equation. Moreover, we provide an iteration method
to find convergent sequences which converge to the solution if one exists.
\end{abstract}
\hrule

\section{Introduction}
Consider the generalized continuous-time algebraic Lyapunov equation (GCALE)
\begin{equation}\label{GCALE}
A^*XB + B^*XA =-Q
\end{equation}
with given matrices $Q, A, B$ and unknown matrix $X$. Such equations play an important role
in stability theory \cite{13, 38}, optimal control problems \cite{33, 37} and balanced model reduction \cite{36}.
In \cite{39}, it  is proved that Eq. (\ref{GCALE}) has a unique Hermitian, positive definite solution $X$ for every Hermitian
positive definite matrix $Q$ if and only if all eigenvalues of the pencil $\lambda A -B$ are finite and lie
in the open left half-plane. When $A=I$ (the identity matrix), Eq. (\ref{GCALE}) becomes the standard  Lyapunov equation
\begin{equation}\label{SLE}
XB + B^*X =-Q.
\end{equation}
The classical numerical methods  to solve   Eq. (\ref{SLE})  are the Bartels-Stewart method \cite{2}, the Hammarling method \cite{19} and the Hessenberg-Schur method \cite{18}. An extension of these methods to solve Eq. (\ref{GCALE}) with  the assumption $A$ is nonsingular, is given in \cite{9, 14, 15, 18, 39}.
 Other approaches to solve Eq. (\ref{GCALE}) are the sign function method \cite{6, 30, 35}, the ADI method \cite{ADI1, ADI2, ADI3}.

In this paper, we consider Eq. (\ref{GCALE}), where $Q$ is an $N\times N$  Hermitian positive definite matrix and $A,B$ are arbitrary $N\times N$ matrices.
Using Bhaskar-Lakshmikantham coupled fixed point theorem \cite{Bhaskar2006},  we provide a sufficient condition that assures the existence and uniqueness of a Hermitian positive definite solution to Eq. (\ref{GCALE}). Moreover, we present an algorithm to solve this equation.  Numerical experiments are given to illustrate our theoretical result.

\section{Notations and preliminaries}
We shall use the following notations: $M(N)$ denotes the set of all $N\times N$  matrices, $H(N)\subset M(N)$ the set of all $N\times N$ Hermitian matrices and
$P(N) \subset H(N)$ is the set of all $N\times N$ positive definite matrices. Instead of $X \in P(N)$, we will also write $X >0$.
Furthermore, $X\geq 0$ means that $X$ is positive semidefinite.
As a different notation for $X-Y \geq 0$ ($X-Y > 0$) we will use $X \geq Y$  ($X > Y$).  If $X, Y \in H(N)$ such
that $X \leq Y$, then $[X,Y]$  will be the set of all $Z \in H(N)$ satisfying $X \leq Z \leq Y$. If $X, Y \in H(N)$ such
that $X < Y$, then $]X,Y[$  will be the set of all $Z \in H(N)$ satisfying $X < Z < Y$.

We denote by $\|\cdot\|$  the spectral norm, i.e., $\|A\|=\sqrt{\lambda^+(A^*A)}$, where $\lambda^+(A^*A)$  is the
largest eigenvalue of $A^*A$. The $N\times N$ identity matrix will be written as $I$.

It turns out that it is convenient here to use the metric induced by the trace norm $\|\cdot\|_1$.
Recall that this norm is given by $\|A\|_1=\sum_{j=1}^n s_j(A)$, where $s_j(A)$, $j = 1,\dots,n$ are the singular values of
$A$. In fact, we shall use a slight modification of this norm. For $Q \in P(N)$,  we define $\|A\|_{1,Q} =\|Q^{1/2}AQ^{1/2}\|_1$.

For any $U\in M(N)$, we denote by $\mbox{Sp}(U)$  the spectrum of $U$, that is, the set of its eigenvalues.

The following lemmas will be useful later.

\begin{lemma}[See \cite{Ran2004}]\label{Lem1}
Let $A \geq 0$ and $B \geq 0$ be $N \times N$ matrices, then $0 \leq tr(AB) \leq \left\| A\right\| \cdot tr(B)$.
\end{lemma}

\begin{lemma}[See \cite{lem}]\label{Lem2}
Let $A \in H(n)$ satisfying $-I<A < I$, then $\|A\| < 1$.
\end{lemma}

\begin{definition}
Let $(\Delta,\le)$ be partially ordered set. We say that a mapping $F:\Delta \times \Delta\to \Delta$ has the mixed monotone property if for all $X,Y,J,K \in \Delta$,
 $$
X \le J,  Y \ge K \Longrightarrow F(X,Y) \le F(J,K).
$$
\end{definition}

The proof of our main result is based on  the following fixed point theorems.

\begin{theorem}[Bhaskar and Lakshmikantham \cite{Bhaskar2006}]\label{TBL}
Let $(\Delta,\leq)$  be a partially ordered set and $d$ be a metric on $X$ such that $(\Delta, d)$ is a complete
metric space. Let $F : \Delta \times \Delta \rightarrow \Delta$ be a continuous mapping having the mixed monotone
property on $\Delta$.  Assume that there exists a $\delta \in [0, 1)$ with
$$
d(F(x,y),F(u,v))\leq \frac{\delta}{2}[d(x,u)+d(y,v)]
$$
for all $x\geq u$ and $y\leq v$.
Suppose also that
\begin{itemize}
\item[(i)]  there exist $x_0, y_0 \in \Delta$  such that $x_0\leq F(x_0,y_0)$ and $y_0\geq F(y_0,x_0)$;
\item[(ii)] every pair of elements has either a lower bound or an upper bound, that is, for every $(x,y)\in \Delta\times \Delta$, there exists a $z\in \Delta$ such that
$x\leq z$ and $y\leq z$.
\end{itemize}
Then, there exists a unique $\overline{x}\in \Delta$ such that $\overline{x}=F(\overline{x},\overline{x})$. Moreover, the sequences $\{x_n\}$ and $\{y_n\}$ defined by $x_{n+1}=F(x_n,y_n)$ and $y_{n+1}=F(y_n,x_n)$ converge to $\overline{x}$, with the following estimate
$$
\max\{d(x_{n},\overline{x}),d(y_{n},\overline{x})\}\leq \frac{\delta^n}{1-\delta} \max\{d(x_0,x_1),d(y_0,y_1)\},\, \mbox{ for all }\,n.
$$
\end{theorem}

\begin{theorem}[Schauder Fixed point theorem]\label{SC}
Let $S$ be a nonempty, compact, convex subset of a normed vector space. Every continuous function $f : S \rightarrow S$ mapping $S$ into itself has a fixed point.
\end{theorem}

\section{Main result}

Our main result is the following.

\begin{theorem}\label{T1}
Suppose that there exists $(\widetilde{Q},M)\in P(N)\times P(N)$ such that
\begin{enumerate}[{\rm (a)}]
\item $2(U\widetilde{Q}U^*+B\widetilde{Q}B^*)<\widetilde{Q}$;
\item $2(V\widetilde{Q}V^*+W\widetilde{Q}W^*)<\widetilde{Q}$;
\item $U^*MU+B^*MB<M-(V^*MV+W^*MW)$;
\item $Q\in ]U^*MU+B^*MB, M-(V^*MV+W^*MW)[$,
\end{enumerate}
where
$$
U=\frac{A-B+I}{\sqrt 2},\quad V=\frac{A+B+I}{\sqrt 2}\quad\mbox{and} \quad W=B-I.
$$
Then,
\begin{enumerate}[{\rm (i)}]
\item Eq. (\ref{GCALE}) has one and only one solution $\widehat{X}\in P(N)$.
\item $\widehat{X}\in [Q-(U^*MU+B^*MB),Q+(V^*MV+W^*MW)]$.
\item Let $(X_n)$ and $(Y_n)$ the sequences defined by  $X_0= 0$, $Y_0=M$, and
\begin{eqnarray}\label{seqxy}
\left\{\begin{array}{lll}
X_{n+1}&=&Q +(V^*X_nV+W^*X_nW)-(U^*Y_nU+B^*Y_nB)\\
Y_{n+1}&=&Q +(V^*Y_nV+W^*Y_nW)-(U^*X_nU+B^*X_nB)
\end{array}.\right.
\end{eqnarray}
We have
\begin{equation}\label{lim}
\lim_{n\rightarrow \infty}\|X_n-\widehat{X}\|=\lim_{n\rightarrow \infty}\|Y_n-\widehat{X}\|=0,
\end{equation}
and the error estimation is given by
\begin{equation}\label{Eerr}
\max\left\{\|X_{n}-\widehat{X}\|,\|Y_{n}-\widehat{X}\|\right\}\leq \frac{\delta^n}{1-\delta}\max\left\{\|X_1-X_0\|,\|Y_1-Y_0\|\right\},
\end{equation}
for all $n$, where $0<\delta<1$.
\end{enumerate}
\end{theorem}

\noindent{\bf Proof.}
It is easy to show that Eq. (\ref{GCALE}) is equivalent to
\begin{equation}\label{GCALE2}
X=Q +(V^*XV+W^*XW)-(U^*XU+B^*XB).
\end{equation}
Consider the continuous mapping $F: H(N)\times H(N)\rightarrow H(N)$ defined by
\begin{equation}\label{F}
F(X,Y)=Q +(V^*XV+W^*XW)-(U^*YU+B^*YB),\quad \mbox{for all}\quad X,Y\in H(N).
\end{equation}
Clearly, Eq. (\ref{GCALE2}) is equivalent to
\begin{equation}\label{GCALE3}
X=F(X,X).
\end{equation}

Let  $X,Y,J,K\in H(N)$ such that $X\leq J$ and $Y\geq K$. Then,
\begin{eqnarray*}
F(X,Y)&=&Q +(V^*XV+W^*XW)-(U^*YU+B^*YB)\\
&\leq & Q +(V^*JV+W^*JW)-(U^*KU+B^*KB)\\
&=&F(J,K).
\end{eqnarray*}
This implies that $F$ is a mixed monotone mapping.

Let $X,Y,J,K \in H(N)$ such that $X \geq J$ and $Y \leq K$.  We have
\begin{align*}
&\|F(X,Y)-F(J,K)\|_{1,\Q}= \tr\left( \Q^{1/2}(F(X,Y)-F(J,K))\Q^{1/2} \right)\\
&=\tr\left( \Q^{1/2}(V^*(X-J)V+W^*(X-J)W+U^*(K-Y)U+B^*(K-Y)B)\Q^{1/2} \right)\\
&=\tr\left( \Q^{1/2}(V^*(X-J)V + W^*(X-J)W)\Q^{1/2} \right)+\tr\left( \Q^{1/2}(U^*(K-Y)U+B^*(K-Y)B)\Q^{1/2} \right).
\end{align*}
On the other hand, using Lemma \ref{Lem1}, we have
\begin{align*}
&\tr\left( \Q^{1/2}(V^*(X-J)V + W^*(X-J)W)\Q^{1/2} \right) =
\tr\left(V\Q V^*(X-J)+W\Q W^*(X-J) \right)\\
&=\tr\left(V\Q V^*(X-J)\Q ^{1/2}\Q ^{-1/2}+W\Q W^*(X-J)\Q ^{1/2}\Q ^{-1/2}\right)\\
&=\tr\left(\Q ^{-1/2}V\Q V^*(X-J)\Q ^{1/2}+\Q ^{-1/2}W\Q W^*(X-J)\Q ^{1/2}\right)\\
&=\tr\left(\Q ^{-1/2}V\Q V^*\Q ^{-1/2}\Q ^{1/2}(X-J)\Q ^{1/2}+\Q ^{-1/2}W\Q W^*\Q ^{-1/2}\Q ^{1/2}(X-J)\Q ^{1/2}\right)\\
&\le \left\|\Q ^{-1/2}V\Q V^*\Q ^{-1/2}+\Q ^{-1/2}W\Q W^*\Q ^{-1/2}\right\| \tr\left(\Q ^{1/2}(X-J)\Q ^{1/2}\right) \\
&= \left\|\Q ^{-1/2}(V\Q V^*+W\Q W^*)\Q ^{-1/2}\right\| \tr\left(\Q ^{1/2}(X-J)\Q ^{1/2}\right) \\
&=\left\|\Q ^{-1/2}(V\Q V^*+W\Q W^*)\Q ^{-1/2}\right\|\|X-J\|_{1,\Q }.
\end{align*}
Thus, we have
\begin{equation}\label{in1}
\tr\left( \Q^{1/2}(V^*(X-J)V + W^*(X-J)W)\Q^{1/2} \right) \leq \left\|\Q^{-1/2}(V\Q V^*+W\Q W^*)\Q^{-1/2}\right\|\|X-J\|_{1,\Q}.
\end{equation}
Similarly, we have
\begin{equation}\label{in2}
\tr\left( \Q^{1/2}(U^*(K-Y)U+B^*(K-Y)B)\Q^{1/2} \right) \leq \left\|\Q^{-1/2}(U\Q U^*+B\Q B^*)\Q^{-1/2}\right\|\|K-Y\|_{1,\Q}.
\end{equation}
Now, using (\ref{in1}) and (\ref{in2}), we get
$$
\|F(X,Y)-F(J,K)\|_{1,\Q  }\leq \frac{\delta }{2} \left(\|X-J\|_{1,\Q}+\|K-Y\|_{1,\Q}\right),
$$
where
$$
\delta=
2\max\left\{\left\|\Q^{-1/2}(V\Q V^*+W\Q W^*)\Q^{-1/2}\right\|,\left\|\Q^{-1/2}(U\Q U^*+B\Q B^*)\Q^{-1/2}\right\|\right\}.
$$
Now, from Lemma \ref{Lem2},   (a) and (b), we have $\delta\in (0,1)$.

Taking $X_0=0$ and $Y_0=M$, from (c) and (d), it follows that
\begin{equation}\label{MM}
X_0< F(X_0,Y_0)\quad\mbox{and}\quad Y_0> F(Y_0,X_0).
\end{equation}

Since all the hypotheses of Theorem \ref{TBL} are satisfied, we deduce that there exists a unique $\widehat{X}\in H(N)$ solution to Eq. (\ref{GCALE3}).
This implies that Eq. (\ref{GCALE}) has a unique solution $\widehat{X}\in H(N)$.

Now, to establish (i), we need to prove that $\widehat{X}\in P(N)$. The Schauder fixed point theorem will be useful in this step.
Define the mapping $G: [F(0,M),F(M,0)]\rightarrow H(N)$ by
$$
G(X)=F(X,X),\quad\mbox{for all}\quad X\in  [F(0,M),F(M,0)].
$$
Note that from (\ref{MM}) and the mixed monotone property of $F$, we have $F(0,M)\leq F(M,0)$.

We shall prove that  $G( [F(0,M),F(M,0)])\subseteq  [F(0,M),F(M,0)]$. Let $X\in [F(0,M),F(M,0)]$, that is,
$$
F(0,M)\leq X\leq F(M,0).
$$
Using the mixed monotone property of $F$, we get
$$
F(F(0,M),F(M,0))\leq F(X,X)=G(X) \leq F(F(M,0),F(0,M)).
$$
On the other hand, from (\ref{MM}), we have
$$
0<F(0,M)\quad\mbox{and}\quad M> F(M,0).
$$
Again, using the mixed monotone property of $F$, we get
$$
F(F(M,0),F(0,M))\leq F(M,0)\quad\mbox{and}\quad F(F(0,M),F(M,0))\geq F(0,M).
$$
Then,
$$
F(0,M) \leq G(X) \leq F(M,0).
$$
Thus we proved that $G( [F(0,M),F(M,0)])\subseteq  [F(0,M),F(M,0)]$.

Now, $G$ maps the compact convex set $[F(0,M),F(M,0)]$ into itself. Since $G$ is
continuous, it follows from Schauder's fixed point theorem (see Theorem \ref{SC})  that $G$ has at least one
fixed point in this set.   However, fixed points of G are solutions of   (\ref{GCALE}), and we
proved already that (\ref{GCALE}) has a unique Hermitian solution. Thus this solution must
be in the set $[F(0,M),F(M,0)]$, that is,
$$
\widehat{X}\in  [Q-(U^*MU+B^*MB),Q+(V^*MV+W^*MW)]\subset P(N).
$$
Thus, we proved (i) and (ii). The proof of (iii) follows immediately from Theorem \ref{TBL}. \hfill    $\square$\\

Now, we present some consequences following from Theorem \ref{T1}, when $A,B$ are Hermitian matrices.

\begin{corollary}\label{CR1}
Suppose that
\begin{enumerate}[{\rm (a)}]
\item $A, B\in H(N)$;
\item $2(UQU+BQB)<Q$;
\item $2(VQV+WQW)<Q$,
\end{enumerate}
where
$$
U=\frac{A-B+I}{\sqrt 2},\quad V=\frac{A+B+I}{\sqrt 2}\quad\mbox{and} \quad W=B-I.
$$
Then,
\begin{enumerate}[{\rm (i)}]
\item Eq. (\ref{GCALE}) has one and only one solution $\widehat{X}\in P(N)$.
\item $\widehat{X}\in [Q-2(UQU+BQB),Q+2(VQV+WQW)]$.
\item Let $(X_n)$ and $(Y_n)$ the sequences defined by  $X_0= 0$, $Y_0=2Q$, and
\begin{eqnarray*}
\left\{\begin{array}{lll}
X_{n+1}&=&Q +(VX_nV+WX_nW)-(UY_nU+BY_nB)\\
Y_{n+1}&=&Q +(VY_nV+WY_nW)-(UX_nU+BX_nB)
\end{array}.\right.
\end{eqnarray*}
We have
\begin{equation*}
\lim_{n\rightarrow \infty}\|X_n-\widehat{X}\|=\lim_{n\rightarrow \infty}\|Y_n-\widehat{X}\|=0,
\end{equation*}
and the error estimation is given by
\begin{equation*}
\max\left\{\|X_{n}-\widehat{X}\|,\|Y_{n}-\widehat{X}\|\right\}\leq \frac{\delta^n}{1-\delta}\max\left\{\|X_1-X_0\|,\|Y_1-Y_0\|\right\},
\end{equation*}
for all $n$, where $0<\delta<1$.
\end{enumerate}
\end{corollary}

\noindent{\bf Proof.}   It follows from Theorem \ref{T1} by taking $\Q = Q$ and $M=2Q$. \hfill   $\square$\\

\begin{corollary}\label{CR2}
Suppose that
\begin{enumerate}[{\rm (a)}]
\item $A,B\in H(N)$;
\item $2(U^2+B^2)<I,\quad 2(V^2+W^2)<I$;
\item $U^2+B^2<Q<I-(V^2+W^2)$,
\end{enumerate}
where
$$
U=\frac{A-B+I}{\sqrt 2},\quad V=\frac{A+B+I}{\sqrt 2}\quad\mbox{and} \quad W=B-I.
$$
Then,
\begin{enumerate}[{\rm (i)}]
\item Eq. (\ref{GCALE}) has one and only one solution $\widehat{X}\in P(N)$.
\item $\widehat{X}\in [Q-(U^2+B^2),Q+(V^2+W^2)]$.
\item Let $(X_n)$ and $(Y_n)$ the sequences defined by  $X_0= 0$, $Y_0=I$, and
\begin{eqnarray*}
\left\{\begin{array}{lll}
X_{n+1}&=&Q +(VX_nV+WX_nW)-(UY_nU+BY_nB)\\
Y_{n+1}&=&Q +(VY_nV+WY_nW)-(UX_nU+BX_nB)
\end{array}.\right.
\end{eqnarray*}
We have
$$
\lim_{n\rightarrow \infty}\|X_n-\widehat{X}\|=\lim_{n\rightarrow \infty}\|Y_n-\widehat{X}\|=0,
$$
and the error estimation is given by
$$
\max\left\{\|X_{n}-\widehat{X}\|,\|Y_{n}-\widehat{X}\|\right\}\leq \frac{\delta^n}{1-\delta}\max\left\{\|X_1-X_0\|,\|Y_1-Y_0\|\right\},
$$
for all $n$, where $0<\delta<1$.
\end{enumerate}
\end{corollary}

\noindent{\bf Proof.}   It follows from Theorem \ref{T1} by taking $\Q = M=I$. \hfill   $\square$

\section{Numerical experiments}

In this section, we present some numerical experiments to check the convergence of the proposed algorithm (\ref{seqxy}).
We take $X_0=0$ and $Y_0=M\in P(N)$. For each iteration $i$, we consider the residual errors
$$
E_i(X)=\|X_i-(Q +(V^*X_iV+W^*X_iW)-(U^*X_iU+B^*X_iB))\|,
$$
$$
E_i(Y)=\|Y_i-(Q +(V^*Y_iV+W^*Y_iW)-(U^*Y_iU+B^*Y_iB))\|
$$
and
$$
E_i=\max\{E_i(X),E_i(Y)\}.
$$
All programs are written in MATLAB version 7.1.

\begin{example}
We consider Eq. (\ref{GCALE}) with
\begin{eqnarray*}
Q=\left(
  \begin{array}{ccc}
  2 &  0.02 &  0.05 \\
  0.02 &     2 &  0.02 \\
  0.05 &  0.02 &     2 \\
  \end{array}
\right),\quad
A=
\left(
  \begin{array}{ccc}
   -0.95 &  0.001 &  0.001 \\
  0.001 & -0.95 &  0.001 \\
  0.001 &  0.001 & -0.95 \\
  \end{array}
\right),\quad
B=\left(
    \begin{array}{ccc}
     0.54 & -0.002 & -0.002 \\
 -0.002 &  0.54 & -0.002 \\
 -0.002 & -0.002 &  0.54 \\
    \end{array}
  \right).
  \end{eqnarray*}
  In this case, we have $A,B\in H(3)$,
$$
\mbox{Sp}\bigg(Q -2(UQU+BQB)\bigg)=\bigg\{0.3345, 0.3379, 0.3880\bigg\},\,\,\mbox{Sp}\bigg(Q -2(VQV+WQW)\bigg)=\bigg\{0.4532, 0.4573, 0.4612\bigg\},
$$
which imply that conditions (a)-(c) of Corollary \ref{CR1} are satisfied.

Considering the iterative method (\ref{seqxy}) with $X_0= 0$ and $Y_0=2Q$, after 100 iterations one gets an approximation to the positive definite solution $\widehat{X}$ and it is

\begin{eqnarray*}
\widehat{X}\approx X_{100}=Y_{100}=\left(
  \begin{array}{ccc}
  1.9495 &  0.0288 &  0.0142 \\
   0.0288 &     1.9496 &  0.0288 \\
  0.0142 &  0.0288 &   1.9495 \\
  \end{array}
\right)
\end{eqnarray*}
and $E_{100}=1.9215\times 10^{-13}$.

\end{example}

\begin{example}
We consider Eq. (\ref{GCALE}) with
\begin{eqnarray*}
Q=\left(
    \begin{array}{ccccc}
      0.4 &    0.01 &    0.02 &    0.03 &    0.04  \\
    0.01 &      0.4  &    0.01 &    0.02 &    0.03  \\
    0.02 &    0.01 &      0.4  &    0.01 &    0.02  \\
    0.03 &    0.02 &    0.01 &      0.4  &    0.01  \\
    0.04 &    0.03 &    0.02 &    0.01 &      0.4 
    \end{array}
  \right),\quad
A=
\left(
    \begin{array}{ccccc}
      -0.95 &   0.001 &   0.001 &   0.001 &   0.001  \\
   0.001 &   -0.95 &   0.001 &   0.001 &   0.001  \\
   0.001 &   0.001 &   -0.95 &   0.001 &   0.001  \\
   0.001 &   0.001 &   0.001 &   -0.95 &   0.001  \\
   0.001 &   0.001 &   0.001 &   0.001 &   -0.95  
    \end{array}
  \right)
  \end{eqnarray*}
and
\begin{eqnarray*}
B=
\left(
    \begin{array}{ccccc}
  0.44 &   -0.02 &   -0.02 &   -0.02 &   -0.02  \\
   -0.02 &    0.44 &   -0.02 &   -0.02 &   -0.02  \\
   -0.02 &   -0.02 &    0.44 &   -0.02 &   -0.02  \\
   -0.02 &   -0.02 &   -0.02 &    0.44 &   -0.02  \\
   -0.02 &   -0.02 &   -0.02 &   -0.02 &    0.44  
    \end{array}
  \right).
\end{eqnarray*}
In this case, we have $A,B\in H(5)$,
$$
\mbox{Sp}\bigg(I -2(U^2+B^2)\bigg)=\bigg\{0.4078,  0.4078, 0.4078,    0.4078,  0.64716\bigg\},
$$
$$
\mbox{Sp}\bigg(I -2(V^2+W^2)\bigg)=\bigg\{0.0094, 0.1577, 0.1577, 0.1577, 0.1577\bigg\},
$$
$$
\mbox{Sp}\bigg(Q -(U^2+B^2)\bigg)=\bigg\{0.0516, 0.0882, 0.0963, 0.0984, 0.3049\bigg\},
$$
and 
$$
\mbox{Sp}\bigg(I -(V^2+W^2)-Q\bigg)=\bigg\{0.0231, 0.1844, 0.1865, 0.1949, 0.2312\bigg\},
$$
which imply that conditions (a)-(c) of Corollary \ref{CR2} are satisfied.

Considering the iterative method (\ref{seqxy}) with $X_0= 0$ and $Y_0=I$, after 82 iterations one gets an approximation to the positive definite solution $\widehat{X}$ and it is

\begin{eqnarray*}
\widehat{X}\approx X_{82}=Y_{82}=
\left(
  \begin{array}{ccccc}
0.4895 &  0.0429 &  0.0541 &  0.0658 &  0.0781 \\
  0.0429 &  0.4878 &  0.0418 &  0.0535 &  0.0658 \\
  0.0541 &  0.0418 &  0.4873 &  0.0418 &  0.0541 \\
  0.0658 &  0.0535 &  0.0418 &  0.4878 &  0.0429 \\
  0.0781 &  0.0658 &  0.0541&  0.0429 &  0.4895 
  \end{array}
\right)
\end{eqnarray*}
and $E_{82}=7.0549\times 10^{-16}$.
\end{example}

\vspace*{1cm}
\noindent Maher Berzig\\
Universit\'e de Tunis,\\
Ecole Sup\'erieure des Sciences et Techniques de Tunis,\\
5, Avenue Taha Hussein-Tunis, B.P. 56, Bab Menara-1008, Tunisie\\
E-mail address: maher.berzig@gmail.com\\

\noindent Bessem Samet\\
Universit\'e de Tunis,\\
Ecole Sup\'erieure des Sciences et Techniques de Tunis,\\
5, Avenue Taha Hussein-Tunis, B.P. 56, Bab Menara-1008, Tunisie\\
E-mail address: bessem.samet@gmail.com


\begin{thebibliography}{99}  {\footnotesize

\bibitem{2}
R.H. Bartels and G.W. Stewart, \textit{ Solution of the equation $AX +XB = C$}, Comm. ACM. 15 (9) (1972) 820-826.



\bibitem{Bhaskar2006}
T.G. Bhaskar and V. Lakshmikantham,
\textit{Fixed point theory in partially ordered metric spaces and applications}, Nonlinear Anal. 65 (2006) 1379--1393.



\bibitem{ADI1}
 P. Benner, J.R. Li and T. Penzl, \textit{ Numerical solution of large Lyapunov equations, Riccati equations, and linear-quadratic control problems}, Numer. Linear Algebra
Appl. 15(9) (2008) 755--777.


\bibitem{6}
P. Benner and E.S. Quintana-Ort\'i, \textit{ Solving stable generalized Lyapunov equations with the matrix sign function},  Numerical Algorithms.  20 (1) (1999) 75--100.


\bibitem{9}
K.E. Chu, \textit{The solution of the matrix equations $AXB-CXD = E$ and $(YA-DZ, YC-BZ)=(E,F)$}, Linear Algebra Appl. 93 (1987) 93--105.


\bibitem{13}
F.R. Gantmacher, \textit{Theory of Matrices}, Chelsea, New York, 1959.

\bibitem{14}
J.D. Gardiner, A.J. Laub, J.J. Amato and C.B. Moler, \textit{ Solution of the Sylvester matrix
equation $AXB^T + CXD^T = E$}, ACM Trans. Math. Software. 18 (2) (1992) 223--231.

\bibitem{15}
J.D. Gardiner, M.R. Wette, A.J. Laub, J.J. Amato  and C.B. Moler, \textit{Algorithm 705: A
Fortran-77 software package for solving the Sylvester matrix equation $AXB^T +CXD^T =E$},  ACM Trans. Math. Software. 18 (2) (1992) 232--238.


\bibitem{18}
G.H.  Golub, S. Nash, and C. Van Loan,  \textit{A Hessenberg-Schur method for the problem $AX + XB = C$}, IEEE Trans. Automat. Control. 24 (1979) 909--913.

\bibitem{19}
S.J. Hammarling, \textit{Numerical solution of the stable non-negative definite Lyapunov equation}, IMA J. Numer. Anal. 2 (1982) 303--323.


\bibitem{ADI2}
J.R. Li and J. White, \textit{Low rank solution of Lyapunov equations} SIAM J. Matrix Anal. Appl. 24(1) (2002) 260--280.


\bibitem{lem}
J.H. Long, X.Y. Hu and L. Zhang, \textit{On the Hermitian positive definite solution
of the nonlinear matrix equation $X + A^*X^{-1}A + B^*X^{-1}B= I$}, Bull Braz Math Soc. New Series.  39(3) (2008)  371--386.


\bibitem{30}
C.S. Kenney and A.J. Laub, \textit{The matrix sign function},  IEEE Trans. Automat. Control. 40 (8) (1995) 1330--1348.


\bibitem{33}
P. Lancaster and L. Rodman, \textit{The Algebraic Riccati Equation}, Oxford University Press, Oxford, 1995.

\bibitem{35}
V.B. Larin and F.A. Aliev, \textit{Generalized Lyapunov equation and factorization of matrix polynomials},  Systems Control Lett. 21 (6) (1993) 485--491.

\bibitem{36}
A.J. Laub, M.T. Heath, C.C. Paige and R.C. Ward, \textit{Computation of system balancing
transformations and other applications of simultaneous diagonalization algorithms},  IEEE Trans. Automat. Control. 32 (2) (1987) 115--122.

\bibitem{37}
V. Mehrmann,  \textit{The Autonomous Linear Quadratic Control Problem, Theory and Numerical Solution}, Lecture Notes in Control and Information Sciences, 163. Springer-Verlag, Heidelberg, 1991.


\bibitem{38}
P.C. M\"{u}ller,  \textit{Stability of linear mechanical systems with holonomic constraints}, Appl. Mech. Rev. 46 (11) (1993) 160--164.

\bibitem{39}
T. Penzl,  \textit{Numerical solution of generalized Lyapunov equations}, Adv. Comput. Math. 8 (1-2) (1998) 33--48.

\bibitem{ADI3}
T. Penzl,  \textit{A cyclic low-rank Smith method for large sparse Lyapunov equations}, SIAM J. Sci. Comp. 21(4)  (1999) 1401--1418.


\bibitem{Ran2004}
A.C.M. Ran and M.C.B. Reurings, \textit{A fixed point theorem in partially ordered sets and some applications to matrix equations}, Proc. Amer. Math. Soc. 132 (2004) 1435--1443.

}


\end{thebibliography}
\end{document}